\documentclass[11pt, letter]{article}

\usepackage{amsmath, amsfonts, amsthm}
\newtheorem{theorem}{Theorem}

\newtheorem{lemma}[theorem]{Lemma}

\def\calV{\mathcal{V}}
\def\calB{\mathcal{B}}
\def\calA{\mathcal{A}}
\def\calC{\mathcal{C}}
\def\ebar{\bar e}
\def\fbar{\bar f}

\def\defeq{:=}

\begin{document}

\title{\bf Two remarks\\ concerning balanced matroids}

\author{Mark Jerrum\thanks{Postal address:  
School of Informatics, University of Edinburgh, 
The King's Buildings, Edinburgh EH9~3JZ, United Kingdom.}\\
School of Informatics\\
University of Edinburgh
}

\maketitle
 
\begin{abstract}\noindent
The property of balance (in the sense of Feder and Mihail) 
is investigated in the context of paving matroids.  
The following examples are exhibited: 
(a)~a class of ``sparse'' paving matroids that are balanced,
but at the same time rich enough combinatorially to
permit the encoding of hard counting problems; and 
(b)~a paving matroid that is not balanced. 
The computational significance of~(a) 
is the following.  As a consequence of balance, there is an 
efficient algorithm for approximating the number of bases 
of a sparse paving matroid within specified relative error. 
On the other hand, determining the number of bases exactly
is likely to be computationally intractable.
\end{abstract}

\section{Discussion}
Let $E$ be a finite ground set and $\calB\subseteq2^E$ a collection of  
subsets of~$E$.  We say that $\calB$ forms the collection of {\it bases\/}
of a {\it matroid} $M=(E,\calB)$ if the following two conditions hold:
\begin{enumerate}
\item All bases (sets in $\calB$) 
have the same size~$r$, namely the {\it rank\/} of $M$.
\item For every pair of bases $X,Y\in\calB$ and every element $e\in X$,
there exists an element $f\in Y$ such that $X\cup\{f\}\setminus\{e\}\in\calB$.
\end{enumerate}
Several other equivalent axiomatisations of matroid are possible,
but the above choice turns out to be the most appropriate for our needs.
For other possible axiomatisations, and more on matroid theory generally,
consult Oxley~\cite{Oxley} or Welsh~\cite{Welsh}. 

The above axioms for a matroid capture the notion of linear independence.
Thus if $S=\{u_0,\ldots,u_{m-1}\}$ is a set of $n$-vectors over a field~$K$,
then the maximal linearly independent subsets of~$S$
form the bases of a matroid with ground set~$S$.  The bases in this instance
have size equal to the dimension of the vector space spanned by~$S$, 
and they clearly satisfy the second or ``exchange'' axiom.
A matroid that arises in this way is {\it vectorial},
and is said to be {\it representable over~$K$}.  
A matroid that is representable over every field is said to
be {\it regular}.  An important combinatorial
example is provided by the collection of 
all spanning trees in an undirected graph $H$: these form the bases of a 
matroid, the {\it cycle matroid of $G$}, with ground set~$E(H)$,
the edge set of~$H$.
A matroid that arises as the cycle matroid of some graph 
is called {\it graphic}.  Graphic matroids form a proper subclass of regular 
matroids.  

The matroid axioms given above suggest a 
natural walk on the set of bases of a matroid~$M$.
The {\it bases-exchange graph\/} $G(M)$ 
of a matroid~$M$ has vertex set~$\calB=\calB(M)$
and edge set 
$$
\big\{\{X,Y\}:X,Y\in\calB\hbox{ and }|X\oplus Y|=2\big\},
$$
where $\oplus$ denotes symmetric difference.  Note that the edges of
the bases-exchange graph~$G(M)$ correspond to the transformations guaranteed
by the exchange axiom.  Indeed, it is straightforward to check,
using the exchange axiom, that the graph $G(M)$ is always connected.
By simulating a random walk on $G(M)$ it is possible, in principle,
to sample a basis almost uniformly at random (u.a.r.\null) from~$\calB(M)$. 
We'll return to this idea presently. 

An intriguing feature of the bases-exchange graph is that it 
appears to have very high ``edge expansion''.  For any matroid $M=(E,\calB)$,
define the {\it edge expansion\/} of its bases-exchange graph to be 
$$
\alpha=\alpha(M)\defeq \min\big\{|\mathrm{cut}(\calA)|/|\calA|:
  \emptyset\subset\calA\subset\calB\text{ and }|\calA|\leq\tfrac12|\calB|\big\},
$$
where $\mathrm{cut}(\calA)$ denotes the cut defined by~$\calA$, 
i.e., the set of edges in~$E$ with one endpoint in~$\calA$   
and one in $\calB\setminus\calA$.  Whenever it has been 
possible to compute the edge expansion of 
the bases-exchange graph of a  matroid~$M$, it 
has been found that~$\alpha(M)\geq1$.  The conjecture that
$\alpha(M)\geq1$ for all matroids is a special case of an 
even stronger one, called the ``zero-one polytope 
conjecture'' of Mihail and Vazirani~\cite{Mihail, Kaibel}.
The circumstantial evidence in favour of the conjecture, even 
in its restricted matroid version, is far from overwhelming. 
Our ignorance concerning the edge expansion of matroids in general 
is almost total:  it is perfectly possible that 
a sequence of matroids exists for which $\alpha(M)$ decays 
exponentially fast as a function of the size~$|E|$ of the ground set.
Nevertheless there is an interesting class of matroids~$M$, the ``balanced''
matroids, for which the lower bound 
$\alpha(M)\geq1$  has been established.  The definition
of balanced matroid (given below) 
is due to Feder and Mihail~\cite{FM}, as is the proof that balance
implies expansion.  

Aside from its intrinsic appeal, the expansion conjecture for 
matroids (and even more so the zero-one polytope conjecture)
has important algorithmic consequences, which arise from the
following considerations. 
Suppose we simulate an unbiased random walk 
on the bases-exchange graph~$G(M)$, with uniform transition  
probabilities (which could be taken as $1/rm$, where $m=|E|$
is the size of the ground set and $r$ the rank).
The walk is ergodic and converges to a stationary distribution
on bases which is uniform.  It is possible, in principle, to use the walk
to sample a basis (almost) u.a.r.\null{} from~$\calB(M)$.  
From there it is a short step (see, e.g., \cite[\S3.2]{ETH}
for the general principle) to estimating the number 
$|\calB(M)|$ of bases within arbitrarily small relative error.
The efficiency of this approach depends crucially on the 
``mixing time'' (number of steps to convergence to near-stationarity)
of the random walk.  As far as we know, this mixing time 
could be exponential in~$m$.  However, the 
mixing time is short whenever $\alpha(M)\geq1$ (or something 
somewhat weaker) holds.  For example, in the case of balanced 
matroids, the mixing time is known to be $O(rm\log r)$:\footnote{To 
avoid trivialities, assume $r\geq2$.}
see Jerrum and Son~\cite{JS}, which improves quantitatively on Feder 
and Mihail~\cite{FM}.  

The standard examples of balanced matroids 
are regular matroids, which were 
shown to be balanced by Feder and Mihail~\cite{FM},
and uniform matroids, which are trivially balanced.  
(The bases of the uniform matroid of rank~$r$ on~$E$
are all $r$-element subsets of~$E$.)
From an algorithmic point of view, this is unfortunate, since the
bases of a regular matroid may be counted exactly 
via linear algebra, and the number of bases of a uniform matroid 
is trivially $\binom{m}{r}$. 
(It can be shown that the bases of a regular matroid are in
1-1 correspondence with the non-singular $r\times r$ submatrices 
of an $r\times m$ unimodular matrix, and that the  number of these
can be computed using the Binet-Cauchy formula.    
Refer to Dyer and Frieze~\cite[\S3.1]{DyerFrieze} for a 
discussion of this topic.)

The first observation in this paper is that paving matroids
from a certain class, which will be called ``sparse'', are all balanced.
(Definitions of ``paving matroid'', and ``sparse paving matroid''
will be given in~\S\ref{sec:balanced}.)
The class of sparse paving matroids is combinatorially 
rich and it is easy to express one's favourite computationally
hard counting problem in terms of counting the bases of a
sparse paving matroid.  This shows that balance is a concept that 
is not entirely devoid of algorithmic interest:  specifically, there exists 
a class of matroids~$M$ for which (a)~exact calculation of $|\calB(M)|$
is as hard as counting satisfying assignments to a Boolean formula,
and hence is almost certainly computationally intractable, whereas
(b)~$|\calB(M)|$ may be approximated within arbitrarily small relative 
error in polynomial time by simulating random walks 
on bases-exchange graphs.\footnote{The issue of balance is not 
in fact crucial here, as Ch{\'a}vez Lomel{\'\i} 
and Welsh~\cite{ChavezWelsh} have presented a polynomial-time algorithm 
for approximately counting bases of an arbitrary paving matroid.}

The second observation resolves an obvious question raised by the first:
namely, are all paving matroids balanced?  It transpires that 
the answer is no, but the construction of a counterexample
requires non-trivial effort.  The counterexample is
based on the Steiner system $S(5,8,24)$.  Welsh~\cite[\S12.6]{Welsh}
has noted the special position that this Steiner system holds in
the theory of matroids.

A closing historical remark.  Dirk Vertigan (personal communication)
has described a class of balanced matroids, unrelated to paving
matroids, whose bases are hard to count in the sense we have
in mind in this note (and which will be clarified in~\S\ref{sec:hard}).
His result was presented during the DIMACS Special Year on Graph Theory 
and Algorithms (1991--2), but was never published.  
Aside from applying to a different class of matroids, 
his construction was apparently more complicated than the
one given here.  So even if the result is not completely 
new, it seems worthwhile to record it here.

\section{A class of balanced paving matroids}\label{sec:balanced}
Suppose $M=(E,\calB)$ is a matroid of rank~$r$.  
A subset of~$E$ is called an {\it independent set\/} 
if it a subset of some basis in~$\calB$. A subset of~$E$ 
that is not an independent set is a {\it dependent set}. 
A minimal (with respect to set inclusion) dependent set
is a {\it circuit}.
The matroid $M$ is said to be {\it paving\/} if all $(r-1)$-element 
subsets of~$E$ are independent sets.  Alternatively, one could say 
that all circuits of~$M$ are of size either $r$ or~$r+1$.
Every $r$-element subset of~$M$ is thus either a basis or a circuit.

A element of $E$ that is contained in no basis of~$M$ is a {\it loop},
and one that is contained in every basis is a {\it coloop}.
Two absolutely central operations on matroids are contraction and deletion.
Assume that $e\in E(M)$ is neither a loop nor a coloop.
If $e$ is an element of the ground set of~$M$ then the 
matroid $M\setminus e$ obtained by {\it deleting\/}~$e$ has 
ground set $E(M\setminus e)=E(M)\setminus\{e\}$ and  
bases $\calB(M\setminus e)=\{X\subseteq E(M\setminus e): X\in\calB(M)\}$;
the matroid $M/e$ obtained by {\it contracting\/}~$e$ has 
ground set $E(M/e)=E(M)\setminus\{e\}$ and  
bases $\calB(M/e)=\{X\subseteq E(M/e): X\cup\{e\}\in\calB(M)\}$.
Any matroid obtained from~$M$ by a series of contractions and deletions
is a {\it minor\/} of~$M$.

The matroid~$M$ is said to possess the {\it negative correlation property\/}
if the inequality $\Pr(e\in X\mid f\in X)\leq\Pr(e\in X)$ holds 
for all pairs of distinct elements $e,f\in X$, 
where we assume that $X\in\calB$ is chosen u.a.r. 
In other words the knowledge that $f$ is 
present in~$X$ makes the presence of~$e$ less likely.\footnote{We assume 
here that $\Pr(f\in X)>0$, i.e., that $f$ is not a loop.}
Further, the matroid~$M$ is said to be {\it balanced\/} if 
all minors of~$M$ (including~$M$ itself) possess the 
negative correlation property.  For more on balanced matroids
in a general matroidal context, 
refer to Choe and Wagner~\cite{ChoeWagner}.

Let $(E,\calB)$ be a paving matroid of rank~$r$ on ground set~$E$.
We have seen that 
such a matroid is defined by the set $\calC_r$ of circuits 
with $r$~elements.
Oxley~\cite[Prop.~1.3.10]{Oxley} provides the following useful
characterisation of paving matroids.

\begin{lemma}\label{lem:PavingCharacter}
Let $\calC_r\subset 2^E$ be a collection of $r$-element subsets
of~$E$.  Then $\calC_r$ defines (in the above sense)
a paving matroid on~$E$ precisely if the following condition holds:
for all $C,C'\in \calC_r$, if\/ $|C\oplus C'|=2$ then every
$r$-element subset of $C\cup C'$ is in~$\calC_r$.
\end{lemma}

We say that a paving matroid is {\it sparse}\footnote{Clearly, 
the qualifier ``sparse'' is intended to refer to the circuits 
and not the bases of the matroid.} if 
\begin{equation}\label{eq:SparseDef}
|C\oplus C'|>2 \text{ for all distinct circuits } C,C'\in\calC_r.
\end{equation}
Note that, by Lemma~\ref{lem:PavingCharacter}, any collection~$\calC_r$ 
of $r$-element subsets of~$E$ 
satisfying~(\ref{eq:SparseDef}) defines
a (sparse) paving matroid.

\begin{lemma}\label{lem:SparseImpliesBalance}
Sparse paving matroids are balanced.
\end{lemma}
\begin{proof}We first verify that every minor of a sparse 
paving matroid is a sparse paving matroid.  This is routine.
Suppose $M=(E,\calC_r)$
is a sparse paving matroid, and $e\in E$ is arbitrary. 
Note that~$e$ cannot be a coloop (except in the trivial case
$r\geq|E|-1$) and so the rank of $M\setminus e$ is~$r$.
The circuits of size~$r$ in $M\setminus e$ are simply all 
the sets in~$\calC_r$ that avoid~$e$~\cite[Eq.~(3.1.14)]{Oxley}.
This subcollection of~$\calC_r$
clearly continues to satisfy~(\ref{eq:SparseDef}).  Furthermore,
$e$ cannot be a loop (except in the trivial case $r=1$), 
and so the rank of $M/e$ is~$r-1$.  The circuits of size 
$r-1$ in $M/e$ are all $r-1$ element subsets 
$C\subseteq E\setminus\{e\}$
satisfying $C\cup\{e\}\in\calC_r$~\cite[Prop.~3.1.11]{Oxley}.
Again, it is clear that 
this collection of circuits satisfies~(\ref{eq:SparseDef}).

It remains to show that the events $e\in X$ and $f\in X$ are
negatively correlated for all distinct $e,f\in E$, assuming 
$X$ is a basis selected u.a.r.  
Partition 
$\calB=\calB_{ef}\cup\calB_{e\fbar}\cup\calB_{\ebar f}\cup\calB_{\ebar\fbar}$ 
into the sets of bases (respectively)
including $e$ and~$f$, including~$e$ but excluding~$f$,
including~$f$ but excluding~$e$, and excluding both $e$ and~$f$.
Consider the bipartite
subgraph of the bases exchange graph $G(M)$
induced by vertex sets $\calB_{ef}$ and $\calB_{e\fbar}$.
Each vertex (basis) $X\in\calB_{ef}$ is adjacent to at least
$m-r-1$ vertices in~$\calB_{e\fbar}$.  
(Consider the collection of $r$-elements sets 
$\{X\cup\{g\}\setminus\{f\}:g\in E\setminus X\}$.  
Condition~(\ref{eq:SparseDef}) ensures that this collection contains 
at most one circuit.)  On the other hand, each vertex (basis) 
in $X\in\calB_{e\fbar}$ is adjacent to at most $r-1$ vertices 
in~$\calB_{ef}$. (The vertices adjacent to~$X$ are all
of the form $X\cup\{f\}\setminus\{g\}$ for some $g\in X\setminus\{e\}$.)
It follows that 
\begin{equation}
  \label{eq:Correspond1}
  (m-r-1)\,|\calB_{ef}|\leq (r-1)\,|\calB_{e\fbar}|.
\end{equation}

Likewise, consider the bipartite
subgraph induced by vertex sets $\calB_{\ebar\fbar}$ and 
$\calB_{\ebar f}$.  Every vertex $X\in\calB_{\ebar\fbar}$
is adjacent to at least $r-1$ vertices in~$\calB_{\ebar f}$.
(Consider the collection of $r$-elements sets 
$\{X\cup\{f\}\setminus\{g\}:g\in X\}$.
As before, this collection contains at most one circuit.)
On the other hand, each vertex (basis) 
in $X\in\calB_{\ebar f}$ is adjacent to at most $m-r-1$ vertices 
in~$\calB_{\ebar\fbar}$.  
It follows that 
\begin{equation}
  \label{eq:Correspond2}
  (r-1)\,|\calB_{\ebar\fbar}|\leq (m-r-1)\,|\calB_{\ebar f}|.
\end{equation}
(Inequality (\ref{eq:Correspond2}) does not rely on sparseness,
and holds in fact for any paving matroid.) 

Multiplying inequalities (\ref{eq:Correspond1}) and~(\ref{eq:Correspond2})
yields 
$$
|\calB_{ef}|\times |\calB_{\ebar\fbar}|\leq 
   |\calB_{e\fbar}|\times|\calB_{\ebar f}|.
$$
A little algebraic manipulation reveals that this inequality is
equivalent to 
$\Pr(e\in X\mid f\in X)\leq \Pr(e\in X)$ 
where $X$ is selected u.a.r.\null{} from the bases of~$M$.
\end{proof}

It is interesting to note that a simple bound on the density 
of bases of a matroid is sufficient to establish the 
negative correlation property.  Specifically, it suffices that
$|\calB|\geq\big(1-\frac{m-r}{2m^2}\big)\binom mr$.
(Martin Dyer, personal communication.)  However, since balance 
requires negative correlation to hold for all minors, 
it is likely that bases need to be somewhat uniformly
distributed as well as dense.  The sparse paving definition
is a convenient way of ensuring these conditions. 

\section{Counting bases is hard,\\ even in balanced matroids}\label{sec:hard}
In discussing algorithms for matroids, the issue of representation
is necessarily problematic, not least because the number of matroids on 
a ground set of size~$m$ is doubly exponential in~$m$.      
Indeed, it is easy to see that the number of sparse 
paving matroids is already doubly exponential.\footnote{Piff 
and Welsh's lower bound on combinatorial geometries~\cite{PiffWelsh} 
is essentially based on counting sparse paving matroids.}
We may note, in passing,
that regular matroids form only a tiny fraction
of all balanced matroids, since the number of regular matroids 
is only singly exponential in~$m$.

In this section, we avoid the issue of representing instances
of paving matroids by not providing a formal definition
of the bases counting problem.  Instead we indicate a simple
method of encoding hard counting problems in sparse paving matroids,
which hopefully will seem quite natural.  A ``hard'' counting 
problem in this context is one that is \#P-complete.  The class~\#P was 
introduced by Valiant as a counting analogue of the more 
familiar class~NP of decision problems.  He showed~\cite{Valiant}
that many natural counting problems are complete for \#P with 
respect to polynomial-time Turing reducibility, and hence
almost certainly computationally intractable.  In particular,
\#P-completeness provides strong evidence against the existence
of a polynomial-time algorithm.

One of the original problems on Valiant's list 
of \#P-complete problems is counting Hamiltonian cycles
in an undirected graph.
Suppose $H=([r],E)$ is an undirected graph on $r$~vertices
with edge set~$E$. Let $\calC_r$ be the collection of all 
Hamilton cycles in~$H$.  Since any pair of Hamilton cycles 
differ in at least four edges, the collection~$\calC_r$
satisfies~(\ref{eq:SparseDef}), and hence defines 
a sparse paving matroid~$(E,\calB)$ of rank~$r$ on~$E$.
Furthermore, it is clear that the number of Hamiltonian cycles in~$H$
is equal to $|\calC_r|=\binom{m}{r}-|\calB|$.
This gives a natural --- in a general combinatorial, though 
not specifically matroidal sense --- encoding of a \#P-complete 
problem as a sparse paving matroid.

It is interesting to observe that the number $|\calB|$ of bases 
of the matroid just constructed can be efficiently approximated 
(by virtue of balance, or by appeal to~\cite{ChavezWelsh}) 
whereas the number~$|\calC_r|$ of non-bases cannot 
(since even deciding emptiness of~$\calC_r$ is hard).

\section{A paving matroid that is not balanced}
Given that a relatively large subset of paving matroids are
balanced, it is natural to ask whether all paving matroids 
are balanced.  The answer is ``no'', but one has to work 
a little to obtain a counterexample.  The problem is to 
squeeze in enough circuits to violate the negative correlation
property.  

We construct a paving matroid of rank six on a ground set~$E$ 
of size~$24$ containing two distinguished elements $e,f$,
and prove that $e$ and~$f$ are positively correlated.
The construction is based on
the Steiner system $S(5,8,24)$.  Denote by~$E$ the ground set 
of $S(5,8,24)$ and by $\calV\subset 2^E$ its set of blocks.
(The ground sets of the Steiner system and 
of the paving matroid will coincide, so it is notationally convenient
to confuse the two.)
The salient features of $S(5,8,24)$ 
are the following~\cite[\S3.6]{BiggsWhite}:
\begin{itemize}
\item $|E|=24$;
\item each block in $\calV$ is of size eight;
\item each subset of five elements of~$E$ is contained in a 
unique block of~$\calV$.
\end{itemize}

We'll define the desired paving matroid $(E,\calB)$ of rank six 
by specifying its circuits of size six.  Let $e,f$ be 
distinguished elements of~$E$.
We'll declare a subset $C\subset E$ of size six
to be a circuit of the matroid if  
there exists a block $V\in \calV$ with $|V\cap\{e,f\}|=1$ 
and $V\supset C$.  Note that two distinct blocks can have at
most four elements in common, and the same is true of circuits
coming from different blocks.
Hence, by ~Lemma~\ref{lem:PavingCharacter}, these 
circuits define a paving matroid.  The bases of this paving matroid
are simply all six-element sets that are not circuits.

Recall the partition 
$\calB=\calB_{ef}\cup\calB_{e\fbar}\cup\calB_{\ebar f}\cup\calB_{\ebar\fbar}$.
We'll calculate the sizes of the various sets occurring in this
partition, and show that 
$|\calB_{ef}|\times|\calB_{\ebar\fbar}|>|\calB_{e\fbar}|\times|\calB_{\ebar f}|$.
It follows directly that the events $e\in X$ and $f\in X$ are
positively correlated, assuming $X\in\calB$ is selected u.a.r. 
\begin{itemize}
\item $|\calB_{ef}|=7315$.  By construction, there are no 
circuits of size six containing both $e$ and~$f$.  
In other words, every six-element set containing $e$ and~$f$ 
is a base, and $|\calB_{ef}|=\binom{22}{4}=7315$.
 
\item $|\calB_{e\fbar}|=|\calB_{\ebar f}|=22638$.
Denote by $\calV_e\subset\calV$ the blocks of the Steiner system
containing~$e$, and by $\calV_{e\fbar}$ the blocks containing~$e$
but excluding~$f$, etc.
First we count the blocks $\calV_{e\fbar}\subset\calV$ 
containing $e$ but not~$f$.
(Every circuit of size six that contains $e$ but not~$f$
must be a subset of a unique such block.) 
Consider any set $A\subset E\setminus\{e\}$ of size four.
Observe that $A\cup\{e\}$ defines a unique block in $\calV_e$,
and moreover that every such block is defined by exactly 
$\binom{7}{4}$ such sets~$A$.  Thus 
$|\calV_e|=\binom{23}{4}\big/\binom{7}{4}=253$.
Similarly each set $A\subset E\setminus\{e,f\}$ of size three 
defines a unique block in $\calV_{ef}$,
and every such block is defined by exactly 
$\binom{6}{3}$ such sets~$A$.
Thus $|\calV_{ef}|=\binom{22}{3}\big/\binom{6}{3}=77$.
Subtracting, $|\calV_{e\fbar}|=176$.
As we observed, each circuit of size six containing~$e$ but not~$f$  
is contained in a unique block in $\calV_{e\fbar}$.
The number of such circuits is thus $176\times\binom{7}{5}=3696$.
Every six-element set containing~$e$ but excluding~$f$
is a basis unless it is one of these 3696 circuits.
Thus $|\calB_{e\fbar}|=\binom{22}{5}-3696=22638$.
Naturally, $|\calB_{\ebar f}|=|\calB_{e\fbar}|$ by symmetry.
\item $|\calB_{\ebar\fbar}|=72149$.
Every six-element set avoiding both $e$ and $f$ is a base
unless it is contained in a block in $\calV_{e\fbar}$ or $\calV_{\ebar f}$.
Thus  
$|\calB_{\ebar\fbar}|=\binom{22}{6} - 2\times176\times\binom{7}{6}=
72149$.
\end{itemize}
In summary,
$$
|\calB_{ef}|\times|\calB_{\ebar\fbar}|=\frac{89015}{86436}\times
|\calB_{e\fbar}|\times|\calB_{\ebar f}|.
$$

This example is perhaps a little larger than necessary, 
but there are limits to how much it can be simplified.
For example, Wagner~\cite{Wagner} 
shows that any matroid that is not balanced 
must have rank at least four.
Furthermore, in order to violate the negative correlation property,
it is necessary to achieve a high density of circuits of size~$r$.

\subsection*{Acknowledgement}I thank Martin Dyer for enlightening 
and encouraging remarks on an early version of this note.
\bibliography{paving}

\begin{thebibliography}{10}

\bibitem{BiggsWhite}
N.~L. Biggs and A.~T. White.
\newblock {\em Permutation groups and combinatorial structures}, volume~33 of
  {\em London Mathematical Society Lecture Note Series}.
\newblock Cambridge University Press, Cambridge, 1979.

\bibitem{ChavezWelsh}
Laura Ch{\'a}vez~Lomel{\'{\i}} and Dominic Welsh.
\newblock Randomised approximation of the number of bases.
\newblock In {\em Matroid theory (Seattle, WA, 1995)}, volume 197 of {\em
  Contemp. Math.}, pages 371--376. Amer. Math. Soc., Providence, RI, 1996.

\bibitem{ChoeWagner}
Young-Bin Choe and David~G. Wagner.
\newblock {Rayleigh matroids}.
\newblock arXiv:math.CO/0307096.

\bibitem{DyerFrieze}
Martin Dyer and Alan Frieze.
\newblock Random walks, totally unimodular matrices, and a randomised dual
  simplex algorithm.
\newblock {\em Math. Programming}, 64(1, Ser. A):1--16, 1994.

\bibitem{FM}
Tom\'as Feder and Milena Mihail.
\newblock Balanced matroids.
\newblock In {\em Proceedings of the 24th Annual ACM Symposium on Theory of
  Computing (STOC)}, pages 26--38. ACM Press, 1992.

\bibitem{ETH}
Mark Jerrum.
\newblock {\em Counting, sampling and integrating: algorithms and complexity}.
\newblock Lectures in Mathematics ETH Z\"urich. Birkh\"auser Verlag, Basel,
  2003.

\bibitem{JS}
Mark Jerrum and Jung-Bae Son.
\newblock Spectral gap and log-{S}obolev constant for balanced matroids.
\newblock In {\em Proceedings of the 43rd IEEE Symposium on Foundations of
  Computer Science (FOCS'02)}, pages 721--729. IEEE Computer Society Press,
  2002.

\bibitem{Kaibel}
Volker Kaibel.
\newblock {On the expansion of graphs of 0/1-polytopes}.
\newblock arXiv:math.CO/0112146.

\bibitem{Mihail}
Milena Mihail.
\newblock On the expansion of combinatorial polytopes.
\newblock In {\em Mathematical foundations of computer science 1992 (Prague,
  1992)}, volume 629 of {\em Lecture Notes in Comput. Sci.}, pages 37--49.
  Springer, Berlin, 1992.

\bibitem{Oxley}
James~G. Oxley.
\newblock {\em Matroid theory}.
\newblock Oxford Science Publications. The Clarendon Press Oxford University
  Press, New York, 1992.

\bibitem{PiffWelsh}
M.~J. Piff and D.~J.~A. Welsh.
\newblock The number of combinatorial geometries.
\newblock {\em Bull. London Math. Soc.}, 3:55--56, 1971.

\bibitem{Valiant}
Leslie~G. Valiant.
\newblock The complexity of enumeration and reliability problems.
\newblock {\em SIAM J. Comput.}, 8(3):410--421, 1979.

\bibitem{Wagner}
David~G. Wagner.
\newblock {Rank three matroids are Rayleigh}.
\newblock arXiv:math.CO/0403216.

\bibitem{Welsh}
D.~J.~A. Welsh.
\newblock {\em Matroid theory}.
\newblock Academic Press [Harcourt Brace Jovanovich Publishers], London, 1976.
\newblock L. M. S. Monographs, No. 8.

\end{thebibliography}
\end{document}